\documentclass[12pt,a4paper]{article}
\usepackage{amssymb,amsthm,amsmath}
\usepackage{fullpage}
\newcommand{\raw}{\rightarrow}

\newcommand{\lr}{\longrightarrow}
\newcommand{\nn}{\nonumber}
\def\bra#1{\left\langle #1\right|}
\def\ket#1{\left| #1\right\rangle}
\def\hs#1#2{\left\langle #1|#2\right\rangle}
\def\hst#1#2{\left\langle#1\overset{.}{\otimes}#2\right\rangle}
\def\hstb#1#2{\left\langle#1\overset{.}{\otimes}_{B}#2\right\rangle}
\def\hstq#1#2{\left\langle#1\overset{.}{\otimes}_{A(S^4_q)}#2\right\rangle}
\def\hsts#1#2{\left\langle#1\overset{.}{\otimes}_{A(S^7_q)}#2\right\rangle}

\DeclareMathOperator{\Mat}{Mat}
\DeclareMathOperator{\Tr}{Tr}
\DeclareMathOperator{\ch}{ch}
\DeclareMathOperator{\tr}{tr}
\DeclareMathOperator{\SU}{SU}

\newcommand{\uno}{x_{\mbox{\tiny{1}}}}
\newcommand{\unob}{\bar{x}^{\mbox{\tiny{1}}}}
\newcommand{\due}{x_{\mbox{\tiny{2}}}}
\newcommand{\dueb}{\bar{x}^{\mbox{\tiny{2}}}}
\newcommand{\tre}{x_{\mbox{\tiny{3}}}}
\newcommand{\treb}{\bar{x}^{\mbox{\tiny{3}}}}
\newcommand{\qu}{x_{\mbox{\tiny{4}}}}
\newcommand{\qub}{\bar{x}^{\mbox{\tiny{4}}}}

\newcommand{\IR}{\ensuremath{\mathbb{R}}}
\newcommand{\IC}{\ensuremath{\mathbb{C}}}
\newcommand{\IH}{\ensuremath{\mathbb{H}}}
\newcommand{\IN}{\ensuremath{\mathbb{N}}}
\newcommand{\II}{\ensuremath{\mathbb{I}}}
\newcommand{\IP}{\ensuremath{\mathbb{P}}}
\newcommand{\ck}{\ensuremath{\mathcal{K}}}
\newcommand{\ce}{\ensuremath{\mathcal{E}}}
\newcommand{\ca}{\ensuremath{\mathcal{A}}}

\newcommand{\cc}{\ensuremath{\mathcal{C}}}
\newcommand{\hil}{\ensuremath{\mathcal{H}}}
\newcommand{\ot}{\otimes}
\newcommand{\ota}{\otimes_{A(S^4_q)}}

\newcommand{\pot}{\overset{.}{\otimes}}
\newcommand{\op}{\oplus}
\newtheorem{prop}{Proposition}
\newtheorem{rema}{Remark}
\newtheorem{defi}{Definition}
\newtheorem{lemma}{Lemma}

\usepackage[english]{babel}

\begin{document}

\title{A Hopf bundle over a quantum four-sphere \\ from the symplectic group}
\date{}

\author{Giovanni Landi$\strut^{1}$, Chiara Pagani$\strut^{2}$, Cesare Reina$^2$
\\[15pt]
$\strut^{1}$ Dipartimento di Matematica e Informatica, Universit\`a 
di Trieste \\
Via A.Valerio 12/1, I-34127 Trieste, Italy \\ and I.N.F.N., Sezione di Napoli,
Napoli, Italy \\ landi@univ.trieste.it
\\[5pt]
$\strut^{2}$ S.I.S.S.A. International School for Advanced Studies, \\ 
Via Beirut
2-4, I-34014 Trieste, Italy \\ pagani@sissa.it, \;\;reina@sissa.it }

\maketitle

\begin{abstract} We construct a quantum version of the $SU(2)$ Hopf bundle
$S^7 \rightarrow S^4$. The quantum sphere
$S^7_q$  arises from the symplectic group $Sp_q(2)$ and a quantum
$4$-sphere $S^4_q$ is obtained via a suitable self-adjoint idempotent $p$ whose
entries generate the algebra
$A(S^4_q)$ of polynomial functions over it. This projection determines a
deformation of an (anti-)instanton bundle over the classical sphere $S^4$. We
compute the fundamental
$K$-homology class of $S^4_q$ and pair it with the class of $p$ in 
the $K$-theory
getting the value $-1$ for the topological charge. There is a right coaction of
$SU_q(2)$ on $S^7_q$ such that the algebra $A(S^7_q)$ is a non trivial quantum principal bundle
over $A(S^4_q)$ with structure quantum group $A(SU_q(2))$.

\vfill
\noindent SISSA-ISAS/50/2004/FM. ESI (2004) 1519.
\\
MSC: 81R60, 16W30, 19D55
\\
Keywords: Noncommutative Geometry, Quantum Groups, Quantum Spheres, Instanton Quantum Bundles.
\end{abstract}
\thispagestyle{empty}
\vfill\eject
\tableofcontents

\vfill\eject
\section{Introduction}

In this paper we study yet another example of how ``quantization removes
degeneracy" by constructing a new quantum version of the Hopf bundle $S^7
\rightarrow S^4$. This is the first outcome of our attempt to generalize to the
quantum case the ADHM construction of $SU(2)$ instantons together
with their moduli
spaces.
 
The $q$-monopole on two dimensional quantum spheres has been constructed in \cite{BM93} 
more than a decade ago. There it was also introduced the general notion of a quantum principal bundle 
with quantum differential calculi, from a geometrical point of view. With universal differential calculi, this notion
was later realised to be equivalent to the one of Hopf-Galois extension (see e.g. \cite{HM}).
An analogous construction for $q$-instantons and their principal bundles has been an open problem ever since.
A step in this direction was taken in \cite{bct} resulting in a bundle which is only a coalgebra extension \cite{bcdt}. 
Here we present a quantum principal instanton bundle which is a honest Hopf-Galois extension. One advantage is that non-universal 
calculi may be constructed on the bundle, as opposite to the case of a coalgebra bundle where there is not such a possibility.

In analogy with the classical case \cite{at}, it is natural
to start with
the quantum version of the (compact) symplectic groups
$A(Sp_q(n))$, i.e. the Hopf algebras generated by matrix elements
$T_i^j$'s with
commutation rules coming  from the
$R$ matrix of
the $C$-series
\cite{frt}. These quantum groups have comodule-subalgebras
$A(S_q^{4n-1})$ yielding
deformations of the algebras of polynomials over the spheres $S^{4n-1}$, which
give more examples of the general construction of quantum homogeneous spaces \cite{BM93}.

The relevant case for us is $n=2$, i.e. the symplectic quantum
$7$-sphere $A(S_q^7)$, which is generated by the matrix elements of
the first and
the last  column of $T$. Indeed, as we will see, $\overline T^4_i \varpropto
T^1_{4-i}$. A similar conjugation occurs for the elements of the
middle columns,
but contrary to what happens at $q=1$, they do not generate a subalgebra. The
algebra $A(S_q^7)$ is the  quantum version of the homogeneous space
$Sp(2)/Sp(1)$
and the injection
$A(S^7_q)\hookrightarrow A(Sp_q(2))$ is a quantum principal bundle with ``structure Hopf algebra''  $A(Sp_q(1))$.

Most importantly, we show that $S_q^7$ is the total space of a
quantum $SU_q(2)$
principal bundle over a quantum $4$-sphere $S_q^4$. 
Unlike the previous construction, this is obviously not a quantum homogeneous structure.
The algebra $A(S_q^4)$  is
constructed as the subalgebra of $A(S_q^7)$ generated by the matrix 
elements of a
self-adjoint projection $p$ which generalizes the anti-instanton of 
charge $-1$.
This projection will be of the form $vv^*$ with $v$ a $4\times 2$  matrix whose
entries are made out of generators of $A(S_q^7)$. The naive 
generalization of the
classical case produces a subalgebra with extra generators which 
vanish at $q=1$.
Luckily enough, there is just one alternative choice of $v$ which 
gives the right
number of generators of an algebra which deforms the algebra of polynomial
functions of $S^4$. At $q=1$ this gives a projection which is gauge 
equivalent to
the standard one.

This good choice becomes even better because there is a natural coaction of
$SU_q(2)$ on $A(S^7_q)$ with coinvariant algebra $A(S^4_q)$ and the injection $A(S^4_q)\hookrightarrow A(S^7_q)$ turns out to be
 a faithfully flat $A(SU_q(2))$-Hopf-Galois extension.

Finally, we set up the stage to compute the charge of our projection and to prove the non triviality of our principal bundle. 
Following a general strategy of noncommutative index theorem \cite{co}, we construct
representations of the algebra $A(S_q^4)$ and the corresponding 
$K$-homology. 
The
analogue of the fundamental class of $S^4$ is given by a non trivial Fredholm
module $\mu$. The natural coupling between $\mu$ and the projection 
$p$ is computed
via the pairing of the corresponding Chern characters $\ch ^*(\mu)\in
HC^*[A(S_q^4)]$ and
$\ch_*(p)\in HC_*[A(S_q^4)]$ in cyclic cohomology and homology respectively
\cite{co}. As expected the result of this pairing, which is an integer by
principle being the index of a Fredholm operator, is actually $-1$ and therefore the bundle is non trivial.

Clearly the example presented in this paper is very special and 
limited, since it
is just a particular anti-instanton of charge $-1$. Indeed our construction is
based on the requirement that the matrix $v$  giving the projection 
is linear in
the generators of $A(S_q^7)$ and such that $v^*v=1$. This is false 
even classically
at generic moduli  and generic  charge, except for the case 
considered here (and
for a similar construction  for the  case of charge $1$). A more 
elaborate strategy
is needed to tackle the general case.

\section{Odd spheres from quantum symplectic groups}

We recall the construction of quantum spheres associated with the 
compact real form
of the quantum symplectic groups
$Sp_q(N, \IC)\; (N=2n)$, the latter being given in
\cite{frt}. Later we shall specialize to the case $N=4$ and the corresponding
$7$-sphere will provide the `total space' of  our quantum Hopf bundle.

\normalsize
\subsection{The quantum groups $Sp_q(N, \IC)$ and $Sp_q(n)$}

The algebra $A(Sp_q(N, \IC))$ is the associative noncommutative 
algebra generated
over the ring of Laurent polynomials $\IC_q:=\IC [q,q^{-1}]$ by the entries
${T_i}^j,\; i,j=1,\dots,N$ of a  matrix $T$ which satisfy RTT equations:
$$ R\;T_1 T_2 = T_2 T_1 R~, \qquad  \quad T_1 =T \ot 1 \; ,
\quad T_2 = 1 \ot T \; .
$$ In components ${(T \ot 1)_{ij}}^{kl} = {T_i}^k {\delta_{j}}^l$. Here the
relevant $N^2 \times N^2$ matrix $R$ is the one for the  $C_N$ series 
and has the
form
\cite{frt},
\begin{eqnarray} R &=& q\sum_{i=1}^N {e_i}^i
\ot {e_i}^i +
\sum_{\stackrel{i,j=1}{i \neq j,{j'}}}^N {e_i}^i \ot {e_j}^j + q^{-1}
\sum_{i=1}^N {e_{i'}}^{i'} \ot {e_i}^i
\nonumber
\\ &&~~~~~~~~~~+ (q -q^{-1}) \sum_{\stackrel{i,j=1}{i >j}}^N {e_i}^j
\ot {e_j}^i -(q -q^{-1})
\sum_{\stackrel{i,j=1}{i >j}}^N q^{\rho_i - \rho_j} \varepsilon_i
\varepsilon_j ~{e_i}^j \ot {e_{i'}}^{j'} \quad ,
\end{eqnarray} where

$i'= N+1 -i$ ~;

${e_i}^j \in M_n(\IC)$  are the elementary  matrices, i.e.
$({e_j}^i)^k_l=\delta_{jl}\delta^{ik}$ ~;

$\varepsilon_i =1, \mbox{ for } i=1,\dots, n$ ~;

$\varepsilon_i = -1 , \mbox{ for }i=n+1 ,\dots, N$ ~;

$(\rho_1,\dots,\rho_N) =(n,n-1,\dots , 1,-1, \dots ,-n)$ ~ .

\bigskip
\noindent The symplectic group structure comes from the matrix
${C_i}^j=q^{\rho_j}\varepsilon_i \delta_{ij'}$ by  imposing  the additional
relations
$$ TCT^tC^{-1}=CT^tC^{-1}T=1 \quad.
$$ The Hopf algebra co-structures $(\Delta, \varepsilon, S)$ of the 
quantum group
$Sp_q(N,\IC)$ are given by
$$
\Delta (T) =T \pot T \; ,  \quad
\varepsilon(T)=I \; , \quad  S(T)=CT^tC^{-1} \; .
$$ In components the antipode explicitly reads
\begin{equation} {S(T)_i}^j = -q^{\rho_{i'}+\rho_j} \varepsilon_i
\varepsilon_{j'}{T_{j'}}^{i'} \; .
\end{equation} At $q=1$ the Hopf algebra  $Sp_q(N,\IC)$ reduces to 
the algebra of
polynomial functions over the symplectic group  $Sp(N,\IC)$.

The compact real form
   $A(Sp_q(n))$ of the quantum group $A(Sp_q(N,\IC))$  is given by 
taking $q\in \IR$
and the anti-involution \cite{frt}
\begin{equation}\label{inv}
\overline{T}=S(T)^t=C^t T (C^{-1})^t \; .
\end{equation}

\subsection{The odd symplectic spheres}

Let us denote
$$ x_i={T_i}^N \;, \quad v^j = {S(T)_N}^j \; , \quad i,j = 1, \dots, N \; .
$$ As we will show, these generators give subalgebras of
$A(Sp_q(N,\IC))$. With the natural involution \eqref{inv}, the 
algebra  generated
by the
$\{x_i,~v^j\}$  can be thought of  as the algebra $A(S^{4n-1}_q)$ of polynomial
functions  on a quantum sphere of `dimension'
$4n-1$.

  From here on, whenever no confusion arises, the sum over repeated indexes is
understood. In components the  RTT equations are given by
\begin{equation}\label{rtt} {R_{ij}}^{kp} ~{T_k}^r ~{T_p}^s = {T_j}^p ~{T_i}^m
~{R_{mp}}^{rs} \; .
\end{equation} Hence
$$ {R_{ij}}^{kl} ~{T_k}^r = {T_j}^p~{T_i}^m ~{R_{mp}}^{rs} {S(T)_s}^l \; ,
$$ and in turn
$$
\begin{array}{l} {S(T)_p}^j~{R_{ij}}^{kl} = {T_i}^a~{R_{ap}}^{rs}
~{S(T)_s}^l~{S(T)_r}^k \; ,
\end{array}
$$ so that
\begin{equation}\label{ss} {S(T)_a}^i~{S(T)_p}^j~{R_{ij}}^{kl} = {R_{ap}}^{rs}
~{S(T)_s}^l~{S(T)_r}^k \; .
\end{equation} Conversely, if we multiply
${R_{ij}}^{kp}~{T_k}^r = {T_j}^l~{T_i}^m ~{R_{ml}}^{rs} ~ {S(T)_s}^p
$ on the left by $S(T)$ we have
\begin{equation}\label{st} {S(T)_l}^j ~{R_{ij}}^{kp}~{T_k}^r ={T_i}^m 
~{R_{ml}}^{rs}
~{S(T)_s}^p \; .
\end{equation} We shall use equations \eqref{rtt}, \eqref{ss} and \eqref{st} to
describe the algebra generated by the $x_i$'s and by the $v^i$'s.

\subsubsection*{The algebra  $\IC_q[x_i]$}

   From \eqref{rtt} with $r=s=N$ we have
\begin{equation} {R_{ij}}^{kp}~ x_k x_p = {T_j}^p ~{T_i}^m ~{R_{mp}}^{NN} \;.
\end{equation} Since the only element ${R_{mp}}^{NN} ~ \varpropto ~{e_m}^N \ot
{e_p}^N  ~(m,p\leq N)$  which is different from zero is ${R_{NN}}^{NN}=q$,  it
follows that
\begin{equation}\label{comxx}
    {R_{ij}}^{kp}~ x_k x_p = q ~x_j x_i \; ,
\end{equation} and the elements $x_i$'s give an algebra with 
commutation relations
\begin{eqnarray} && x_i x_j = q x_j x_i \; , \quad i<j,~~i\neq j' \; , \nn \\ &&
x_{i'}x_i = q^{-2} x_i x_{i'} + (q^{-2}-1) \sum_{k=1}^{i-1} q^{\rho_i - \rho_k}
\varepsilon_i
\varepsilon_k ~x_k x_{k'} \; , \quad i< i' \; . \label{xx}
\end{eqnarray}

\subsubsection*{The  algebra $\IC_q[v^i]$}\label{v}

Putting  $a=p=N$ in equation \eqref{ss}, we get
$$ v^i v^j {R_{ij}}^{kl}= {R_{NN}}^{rs} {S(T)_s}^l {S(T)_r}^k \; .
$$ The sum on the r.h.s. reduces to ${R_{NN}}^{NN} {S(T)_N}^l 
{S(T)_N}^k$ and the
$v^i$'s give an algebra with commutation relations
\begin{equation}\label{commvv} v^l v^k {R_{lk}}^{ji} = q v^i v^j .
\end{equation} Explicitly
\begin{eqnarray} && v^i v^j = q^{-1} v^j v^i \; , \quad i<j,~~i\neq 
j' \; , \nn \\
&& v^{i'} v^i = q^{2} v^i v^{i'} + (q^{2}-1) \sum_{k=i'+1}^{N} q^{\rho_k -
\rho_{i'}} ~\varepsilon_k \varepsilon_{i'} v^k v^{k'} \; , \quad i<i' \; .
\label{vv}
\end{eqnarray}

\subsubsection*{The algebra  $\IC_q[x_i, v^j]$} 
Finally, for  $l=r=N$ 
the equation
\eqref{st} reads:
    $$ v^j {R_{ij}}^{kp} ~x_k = {T_i}^m {R_{mN}}^{Ns} {S(T)_s}^p.$$ 
Once more, the
only term in $R$ of the form ${e_m}^N \ot {e_N}^s ~(m\leq N)$  is  ${e_N}^N \ot
{e_N}^N$ and therefore
\begin{equation}\label{commxv} v^j {R_{ij}}^{kp} ~x_k = q ~x_i v^p  \; .
\end{equation} Explicitly the mixed commutation rules for the algebra 
$\IC_q[x_i,
v^j]$ read,
\begin{eqnarray} && x_i v^i = v^i x_i + (1 - q^{-2}) \sum_{k=1}^{i-1} 
v^k x_k  +
\underbrace{(1-q^{-2}) q^{\rho_i - \rho_{i'}} v^{i'} x_{i'}}_{if 
~~i>i'} \; , \nn \\
&& x_i v^{i'} = q^{-2} v^{i'} x_i \; , \nn \\ && x_i v^j = q^{-1} v^j x^i \; ,
\quad i \neq j \quad {\rm and} \quad i<j' \nn \\ && x_i v^j = q^{-1} v^j x^i+
(q^{-2}-1) q^{\rho_i - \rho_{j'}} ~\varepsilon_i \varepsilon_{j'} 
v^{i'}x_{j'} \; ,
\quad i \neq j \quad {\rm and} \quad i>j' \;. \label{xv}
\end{eqnarray}

\bigskip
\subsubsection*{The quantum spheres $S^{4n-1}_q$}\label{com}

Let us observe that with the anti-involution \eqref{inv} we have the
identification  $v^i={S(T)_N}^i= \bar{x}^i$. The subalgebra $A(S^{4n-1}_q)$ of
$A(Sp_q(n))$ generated by
$\{x_i, v^i=\bar{x}^i, ~i=1,\dots,2n\}$ is the algebra of polynomial 
functions on a
sphere. Indeed
$$  S(T)T=I~ \Rightarrow  ~\sum {S(T)_N}^i {T_i}^N =\delta_N^N =1 \;
$$  i.e.
\begin{equation}\label{sphere}
\sum_i \bar{x}^i x_i = 1 \; .
\end{equation} Furthermore, the restriction of the comultiplication 
is a natural
left coaction
$$
\Delta_L : A(S^{4n-1}_q) \lr A(Sp_q(n)) \ot A(S^{4n-1}_q)\; .
$$ The fact that $\Delta_L$ is an algebra map then implies that 
$A(S^{4n-1}_q)$  is
a comodule algebra over
$A(Sp_q(n))$.

At $q=1$ this algebra reduces to the algebra of polynomial functions over the
spheres $S^{4n-1}$ as homogeneous  spaces of the symplectic group $Sp(n):\;
S^{4n-1}=Sp(n)/Sp(n-1)$.

\subsection{The symplectic 7-sphere $S^7_q$}\label{s7}

The algebra
$A(S^7_q)$ is generated by the elements $x_i = {T_i}^4$ and
$\bar{x}^i={S(T)_4}^i= q^{2+\rho_i}
\varepsilon_{i'}{T_{i'}}^1,$ for $i=1,\dots ,4$.
   From $S(T)~T = 1$ we have the sphere relation $\sum_{i=1}^4
\bar{x}^i x_i = 1$. Since we shall systematically use them in the following, we
shall explicitly give the commutation relations among the generators. \\
   From \eqref{xx}, the algebra of the $x_i$'s is given by
\begin{equation}
\begin{array}{ll}
\uno \due = q \due \uno \; , & \uno \tre = q \tre \uno \; ,
\\
\due \qu = q \qu \due \; , & \tre \qu = q \qu \tre \; ,
\\
\qu \uno = q^{-2} \uno \qu \; , &
\tre \due = q^{-2} \due \tre +q^{-2}(q^{-1}-q) \uno \qu \; ,
\end{array}
\end{equation} together with their conjugates (given in \eqref{vv}). \\ We have
also the commutation relations between the
$x_i$ and the $\bar{x}^j$ deduced  from \eqref{commxv}:
\begin{equation}
\begin{array}{l}
\begin{array}{ll}
\uno  \unob = \unob \uno \; ,&
\uno \dueb = q^{-1} \dueb \uno \; ,
\\
\uno \treb = q^{-1} \treb \uno \; , &
\uno \qub = q^{-2} \qub \uno \; ,
\end{array}
\\
\\
\begin{array}{l}
\due \dueb = \dueb \due + (1-q^{-2})\unob \uno \; ,\\
\due \treb = q^{-2} \treb \due \; , \\
\due \qub = q^{-1} \qub \due +q^{-1} (q^{-2}-1) \treb \uno \; ,
\end{array}
\\
\\
\begin{array}{ll}
\tre \treb = \treb \tre +(1-q^{-2})[\unob \uno +(1+q^{-2})\dueb \due] \; ,
\\
\tre \qub = q^{-1} \qub \tre +(1-q^{-2})q^{-3} \dueb \uno \; ,
\end{array}
\\
\\
\begin{array}{ll}
\qu \qub = \qub \qu +(1-q^{-2})[(1+q^{-4})\unob \uno + \dueb \due + \treb
\tre] \; ,
\end{array}
\end{array}
\end{equation} again with their conjugates.
\\

Next we show that the algebra $A(S^7_q)$ can be realized as the subalgebra of
$A(Sp_q(2))$ generated by the coinvariants  under the right-coaction of
$A(Sp_q(1)),$ in complete analogy with  the classical homogeneous space
$Sp(2)/Sp(1) \simeq S^7$.

\begin{lemma} The two-sided *-ideal in $A(Sp_q(2))$ generated as
$$ I_q = \{{T_1}^1-1,
{T_4}^4-1,{T_1}^2,{T_1}^3,{T_1}^4,{T_2}^1,{T_2}^4,{T_3}^1,{T_3}^4,{T_4}^1,{T_4}^2,{T_4}^3
\}
$$ with the involution \eqref{inv} is a Hopf ideal.
\end{lemma}
\noindent {\it Proof.} Since  ${S(T)_i}^j \propto{T_{j'}}^{i'}$,
$S(I_q)\subseteq
I_q$ which also proves that $I_q$ is a *-ideal.  One easily  shows that
$\varepsilon (I_q)=0$ and $\Delta (I_q)\subseteq I_q\ot
A(Sp_q(2))+A(Sp_q(2))\ot
I_q$. \qed

\begin{prop} The Hopf algebra $B_q := A(Sp_q(2))/I_q$  is isomorphic to the
coordinate algebra $A(SU_{q^2}(2)) \cong A(Sp_q(1))$.
\end{prop}
\noindent {\it Proof.} Using $\overline{T}=S(T)^t$ and setting
${T_2}^2 = \alpha,~
{T_3}^2=\gamma $, the algebra $B_q$ can be described as the algebra
generated by the
entries of the matrix
\begin{equation}\label{spq1} T'=\left(
\begin{array}{cccc} 1& 0 & 0 &0
\\ 0 & \alpha & -q^2 \bar{\gamma}  & 0
\\ 0 & \gamma &  \bar{\alpha} & 0
\\ 0& 0 & 0 & 1
\end{array}
\right).
\end{equation} The commutation relations deduced from RTT equations \eqref{rtt}
read:
\begin{equation}
\begin{array}{ll}
\alpha \bar\gamma = q^2 \bar\gamma \alpha \quad , \quad &
\alpha \gamma = q^2 \gamma \alpha ~,
\quad
\gamma \bar\gamma= \bar\gamma \gamma ~,
\\
\bar\alpha \alpha + \bar\gamma \gamma =1  \quad ;&
\alpha \bar\alpha + q^4 \gamma \bar\gamma =1 ~.
\end{array}
\end{equation} Hence, as an algebra $B_q$  is isomorphic to the algebra
$A(SU_{q^2}(2))$.  Furthermore, the restriction of the coproduct of 
$A(Sp_q(2))$
to $B_q$  endows the latter with a  coalgebra structure, $\Delta(T') 
= T' \pot T'$,
which is the same as the one of
$A(SU_{q^2}(2))$.  We can conclude that also as a Hopf algebra, $B_q$  is
isomorphic to the  Hopf algebra
$A(SU_{q^2}(2)) \cong A(Sp_q(1))$.
\qed
\\
\begin{prop}\label{p2} The algebra  $A(S_q^7) \subset  A(Sp_q(2))$ is
the algebra
of coinvariants with respect to the natural right coaction
$$
\Delta_R:A(Sp_q(2)) \raw A(Sp_q(2)) \pot
A(Sp_{q}(1))
\quad ; \quad
\Delta_R(T) = T \pot T' \; .
$$
\end{prop}
\noindent
{\it Proof.} It is straightforward to show that the generators of the
algebra $A(S_q^7)$ are  coinvariants:
$$
\Delta_R (x_i)= \Delta_R (T_i^4)= x_i \ot 1 ~~;~~
\Delta_R (\bar{x}^i)= -q^{2 + \rho_i} \varepsilon_i ~\Delta_R (T_i^1)=
\bar{x}^i \ot 1 \,
$$
thus the algebra $A(S_q^7)$ is made of coinvariants. There are no other coinvariants of degree one
since each row of the submatrix of $T$ made out of the two central columns is a fundamental comodule under the coaction of $SU_{q^2}(2)$.
Other coinvariants arising at higher even degree are of the form 
$(T_{i2}T_{i3}-q^2 T_{i3}T_{i2})^n$; thanks to the commutation relations of $A(Sp_q(2))$, one checks these belong to $A(S^7_q)$ as well.
It is an easy computation to check that similar expressions involving elements from different rows cannot be coinvariant.
\qed \\

The previous construction is one more example of the general construction of a quantum principal bundle over a quantum homogeneous space \cite{BM93}. 
The latter is the datum of a Hopf quotient $\pi: A(G) \to A(K)$ with the right coaction of $A(K)$ on $A(G)$ given by the reduced coproduct 
$\Delta_R:=(id \ot \pi)\Delta$ where $\Delta$ is the coproduct of $A(G)$. The subalgebra $B \subset A(G)$ made of the coinvariants with respect
to $\Delta_R$ is called a quantum homogeneous space. To prove that it is a quantum principal bundle one needs some more assumptions (see
Lemma 5.2 of \cite{BM93}). In our case $A(G)=A(Sp_q(2)),~A(K)=A(Sp_q(1))$ with $\pi(T)=T'$. We will prove in Sec. \ref{se:hg} that 
the resulting inclusion $ B=A(S^7_q) \hookrightarrow A(Sp_q(2))$ is indeed a Hopf Galois extension and hence a quantum principal bundle.

\section{The principal bundle  $A(S^{4}_q) \hookrightarrow  A(S^{7}_q)$}

The fundamental step of this paper is to make the sphere $S_q^7$
itself into the total space of a quantum principal bundle over a
deformed $4$-sphere. Unlike what we saw in the previous section, this is not a quantum homogeneous
space construction and it is not obvious that such a bundle exists at all.
Nonetheless the notion of quantum bundle is more general and one only needs
that the total space algebra is a comodule algebra over a Hopf algebra with additional suitable properties.

The notion of quantum principle bundle, as said, is encoded in the one of Hopf-Galois extension (see e.g. \cite{BM93}, \cite{HM}).
Let us recall some relevant definitions \cite{KT}  (see also 
\cite{HG}). Recall that we work over the field $k=\IC$.

\begin{defi} Let $H$ be a Hopf algebra and $P$ a right $H$-comodule
algebra with
multiplication $m: P\ot P \raw P$ and  coaction $\Delta_R : P \raw P 
\ot H$. Let
$B\subseteq P$ be the subalgebra of coinvariants, i.e.
$B=\{p\in P ~|~ \Delta_R(p) = p \ot 1\}$. The extension $B\subseteq 
P$ is called an
$H$ Hopf-Galois extension if the canonical map
\begin{eqnarray}\label{can} && \chi : P \otimes_B P \lr P \ot H \; , 
\nn \\ && \chi
:= (m \ot id) \circ (id \otimes_B \Delta_R) \; , \quad p' \ot_B p
\mapsto
\chi(p' \otimes_B p) = p' p_{(0)} \ot p_{(1)} 
\end{eqnarray} is bijective.
\end{defi}
\noindent We use Sweedler-like notation $\Delta_R p=p_{(0)} \ot p_{(1)}$. The 
canonical map
is left $P$-linear and right
$H$-colinear and is a morphism (an isomorphism  for Hopf-Galois 
extensions) of left
$P$-modules and right $H$-comodules. It is also clear that $P$ is 
both a left and a
right $B$-module. \\

The injectivity of the canonical map dualizes the condition of a group  action
$X
\times G
\raw X$ to be free: if $\alpha$ is the map $\alpha: X \times G \raw X
\times_M X, ~(x,g) \mapsto (x,x\cdot g)$ then $\alpha^*=\chi$ with
$P,H$ the algebras of functions on $X,G$ respectively and the action 
is free if and
only if $\alpha$ is injective. Here $M:=X/G$ is the space of orbits 
with projection
map $\pi: X \raw M, ~\pi(x\cdot g)=\pi(x)$, for all
$x\in X, g\in G$. Furthermore, $\alpha$ is surjective if and only if for all
$x\in X$, the fibre $\pi^{-1}(\pi(x))$ of $\pi(x)$ is equal to the 
residue class $x
\cdot G$, that is, if and only if $G$ acts transitively on the fibres of $\pi$.
\\

In differential geometry a principle bundle is more than just a free and effective action
of a Lie group. In our example, thanks to the fact that the ``structure group" is $SU_q(2)$, from Th.~I of \cite{scH} further nice
properties can be established. We shall elaborate more on these points later on in Sect. \ref{se:hg} .
\\

The first natural step would be to construct a map from $S_q^7$ into a
deformation of the
Stieffel variety  of unitary frames of $2$-planes in $\IC ^4$ to parallel the
classical construction as recalled in the Appendix~\ref{se:chf}. The
naive choice
we have is to take as generators the elements of two (conjugate) columns of the
matrix $T$. We are actually forced to take the first and the last
columns of the
matrix $T$ because the other choice (i.e. the second and the third
columns) does
not yield a subalgebra since commutation relations of their elements
will involve
elements from the other two columns. If we set
\begin{equation} v= \left(
\begin{array}{cc}\qub & \uno
\\ q^{-1}\treb & \due
\\
    -q^{-3} \dueb & \tre
\\
    -q^{-4}\unob & \qu
\end{array}
\right) \; ,
\end{equation} we have  $v^*v=\mathbb{I}_2$ and the matrix $p=v~v^*$ is a
self-adjoint idempotent, i.e. $p=p^*=p^2$.  At $q=1$ the entries of $p$ are
invariant for the natural action of $SU(2)$ on $S^7$ and generate the
algebra of
polynomials on $S^4$. This fails to be the case at generic
$q$ due to the occurrence of extra generators e.g.
\begin{equation}
   p_{14}=(1-q^{-2})\uno \qub, \quad p_{23}=(1-q^{-2})\due \treb,
\end{equation} which vanish at $q=1$.

\subsection{The quantum sphere $S_q^4$}

These facts indicate that the naive
quantum analogue of
the quaternionic projective line as a homogeneous space of $Sp_q(2)$
has not the
right number of generators. Rather surprisingly, we shall anyhow be
able to select
another subalgebra of $A(S_q^7)$ which is a deformation of the algebra of
polynomials on $S^4$ having the same number of generators. These
generators come
from a better choice of a projection.

On the free module $\ce:=\IC^4 \ot A(S_q^7)$ we consider the
hermitean structure
given by
$$h(\ket{\xi_1},\ket{\xi_2})=  \sum_{j=1}^4 \bar{\xi_1}^j \xi_2^j \;
.$$ To every
element $\ket{\xi}\in\ce$ one associates an element
$\bra{\xi}$ in the dual module $\ce^*$ by the pairing
$$\bra{\xi}(\ket{\eta}):=\hs{\xi}{\eta}=h(\ket{\xi},\ket{\eta}).$$ Guided 
the classical construction which we present in 
Appendix~\ref{se:chf}, we
shall look for two elements $\ket{\phi_1},~\ket{\phi_2}$  in $\ce$ with the
property that
$$
\hs{\phi_1}{\phi_1} =1 \;, \quad \hs{\phi_2}{\phi_2} =1 \;, \quad
\hs{\phi_1}{\phi_2} =0 \; .
$$ As a consequence, the matrix valued function defined by
\begin{equation} p := \ket{\phi_1} \bra{\phi_1} + \ket{\phi_2}\bra{\phi_2} \; ,
\end{equation} is a self-adjoint idempotent (a projection). \\ In 
principle, $ p
\in \Mat_4(A(S^7_q))$, but we can choose
$\ket{\phi_1},~\ket{\phi_2}$ in such a way that the entries of $p$ 
will generate a
subalgebra $A(S^4_q)$ of
$A(S^7_q)$ which is a deformation of the algebra of polynomial functions on the
$4$-sphere $S^4$. The two elements $\ket{\phi_1},~\ket{\phi_2}$ will 
be obtained in
two steps as follows. \\

Firstly we write the relation $1=\sum \bar{x}^i x_i$  in terms of the quadratic
elements
$\unob \uno$, $\due \dueb$, $\treb \tre$, $\qu \qub$ by using the commutation
relations of Sect.~\ref{s7}. We have that
$$ 1= q^{-6} \unob \uno +q^{-2} \due \dueb + q^{-2} \treb \tre + \qu
\qub \; .$$ Then we take,
\begin{equation}\label{psi1}
\ket{\phi_1}= (q^{-3} \uno ,- q^{-1} \dueb , q^{-1} \tre ,-\qub )^t \; ,
\end{equation} ($t$ denoting transposition) which is such that
$\hs{\phi_1}{\phi_1}=1$. \\

\noindent Next, we write
$1=\sum \bar{x}^i x_i $ as a function of the quadratic elements
$\uno \unob, \dueb \due, \tre \treb, \qub \qu$:
$$ 1= q^{-2} \uno \unob +q^{-4} \dueb \due +  \tre \treb + \qub \qu \; .
$$ By taking,
$$
\ket{\phi_2}= (\pm q^{-2} \due , \pm q^{-1} \unob , \pm  \qu ,\pm \treb )^t
$$ we get $ \hs{\phi_2}{\phi_2}=1$. The signs will be chosen in order 
to have also
the orthogonality $\hs{\phi_1}{\phi_2}=0$; for
\begin{equation}\label{psi2}
\ket{\phi_2}= ( q^{-2} \due ,  q^{-1} \unob ,  - \qu ,- \treb )^t
\end{equation} this is satisfied. \\

\noindent The matrix
\begin{equation}\label{u} v=\left( \ket{\phi_1},  \ket{\phi_2}
\right) = \left(
\begin{array}{cc} q^{-3} \uno & q^{-2} \due
\\ -q^{-1}\dueb & q^{-1} \unob
\\
    q^{-1} \tre & -\qu
\\
    -\qub & -\treb
\end{array}
\right) \; .
\end{equation} is such that $v^*v=1$ and hence $p=vv^*$ is a self-adjoint
projection.

\begin{prop} The entries of the projection
$p= vv^*$, with $v$ given in \eqref{u}, generate a subalgebra of $A(S^7_q)$ 
which is a
deformation of the algebra of polynomial functions on the
$4$-sphere $S^4$.
\end{prop}
\noindent {\it Proof.} Let us compute explicitly the components of 
the projection
$p$ and their commutation relations.
\begin{enumerate}
\item  The diagonal elements are given by
$$
\begin{array}{ll} p_{11}= q^{-6} \uno \unob + q^{-4} \due \dueb \; , 
\; & p_{22}=
q^{-2} \dueb \due  + q^{-2} \unob \uno \; , \\ p_{33} = q^{-2} \tre \treb + \qu
\qub \; , \; & p_{44}= \qub \qu + \treb \tre \; ,
\end{array}
$$ and satisfy the relation
\begin{equation}\label{traccia} q^{-2} p_{11} + q^2 p_{22} + p_{33} + 
p_{44} =2 \; .
\end{equation} Only one of the $p_{ii}$'s is independent; indeed by using the
commutation relations and  the equation $\sum \bar{x}^i x_i=1,$  we 
can rewrite the
$p_{ii}$'s in terms of
\begin{equation} t:= p_{22} \; ,
\end{equation} as
$$ p_{11}=q^{-2}t ~~,~~p_{22}=t ~~,~~p_{33}= 1- q^{-4}t ~~,~~p_{44}= 1-q^2 t
\; .
$$ Equation \eqref{traccia} is easily verified. Notice that $t$ is 
self-adjoint:
$\bar{t}=t$.
\item As in the classical case, the elements $p_{12}, p_{34}$ (and their
conjugates) vanish:
$$ p_{12}=- q^{-4} \uno \due + q^{-3} \due \uno = 0 \; , \quad 
p_{34}=-q^{-1} \tre
\qu + \qu \tre=0 \; .
$$
\item The remaining elements are given by
$$
\begin{array}{ll} p_{13}=q^{-4} \uno \treb - q^{-2} \due \qub \; , 
\qquad & p_{14}=
-q^{-3} \uno \qu -q^{-2} \due \tre  \; , \\ p_{23} =- q^{-2} \dueb 
\treb -q^{-1}
\unob \qub \; , & p_{24} = q^{-1} \dueb \qu - q^{-1} \unob \tre \; ,
\end{array}
$$ with $p_{ji}=\bar{p}_{ij}$ when $j>i$. \\ By using the commutation 
relations of
$A(S^7_q)$, one finds that only two of these are independent. We take 
them to be
$p_{13}$ and $p_{14}$; one finds  $p_{23}= q^{-2}
\bar{p}_{14}$ and $p_{24}=-q^2 \bar{p}_{13}$.
\end{enumerate} Finally, we also have the sphere relation,
\begin{equation}\label{quad} (q^6 - q^8)p_{\mbox{\tiny {11}}}^2+p_{\mbox{\tiny
{22}}}^2 +p_{\mbox{\tiny {44}}}^2+ 
q^4(p_{\mbox{\tiny{13}}}p_{\mbox{\tiny{31}}}+
p_{\mbox{\tiny{14}}}p_{\mbox{\tiny{41}}})+
q^2(p_{\mbox{\tiny{24}}}p_{\mbox{\tiny{42}}}
+p_{\mbox{\tiny{23}}}p_{\mbox{\tiny{32}}}) = (\sum \bar{x}^i x_i)^2 = 1 \; .
\end{equation}

\bigskip\noindent Summing up, together with $t= p_{22}$, we set
$a:=p_{13}$ and $b:=p_{14}$. Then the projection $p$ takes the following form
\begin{equation}\label{proj} p=\left(
\begin{array}{llll} q^{-2}t & 0 & a &b
\\
\\ 0 &t & q^{-2} \bar{b} & -q^2 \bar{a}
\\
\\
\bar{a} & q^{-2} b & 1-q^{-4}t & 0
\\
\\
\bar{b}& -q^2 a & 0 & 1-q^2t
\end{array}
\right) \; .
\end{equation} By construction $p^*=p$ and this means that 
$\bar{t}=t$, as observed,
and that $\bar{a},\bar{b}$ are conjugate to $a,b$ respectively.   Also, by
construction
$p^2=p$;  this property gives the easiest way to compute the 
commutation relations
between the generators. One finds,
\begin{equation}\label{s4}
\begin{array}{ll} ab = q^4 ba \; , & \; \bar{a}b =b \bar{a} \; , \\ 
ta= q^{-2} at
\; , & \;tb = q^4 bt \; ,
\end{array}
\end{equation}
together with their conjugates, and sphere relations
\begin{equation}\label{sr4}
\begin{array}{lll} a\bar{a}+b\bar{b}=q^{-2}t(1-q^{-2}t) \; , ~~~~~ & 
q^4 \bar{a}a
+q^{-4}\bar{b}b=t(1-t) \; ,
\\ b \bar b - q^{-4}\bar{b}b =(1- q^{-4})t^2\; .&
\end{array}
\end{equation} It is straightforward to check also the relation \eqref{quad}.
\qed

\bigskip We define the algebra $A(S^4_q)$ to be  the algebra generated by the
elements $a,\bar{a},b,\bar{b},t$ with the commutation relations \eqref{s4} and
\eqref{sr4}. For $q=1$ it reduces to the algebra of polynomial functions on the
sphere $S^4$. Otherwise, we can limit ourselves to $|q|<1$, because the map
$$q\mapsto q^{-1},\quad a\mapsto q^2\bar a,\quad b\mapsto q^{-2}\bar b,
\quad t\mapsto q^{-2} t$$ yields an isomorphic algebra. \\

At $q=1$, the projection $p$ in \eqref{proj} is conjugate to the classical one
given in Appendix \ref{se:chf} by the matrix $diag[1,-1,1,1]$ (up to a
renaming of the
generators). \\

Our sphere $S^4_q$ seems to be different from the one constructed in
\cite{bct}. Two of our generators commute and most importantly, it 
does not come
from a deformation of a subgroup (let alone coisotropic) of $Sp(2)$. 
However, at
the continuous level these two quantum spheres are the same since  the
$C^*$-algebra completion of both polynomial algebras is the minimal unitization
$\ck \op \IC \II $ of the compact operators on an infinite 
dimensional separable
Hilbert space, a property shared with Podle{\'s} standard sphere as 
well \cite{po}.
This fact will be derived in Sect.~\ref{se:rep} when we study the 
representations
of  the algebra $A(S^4_q)$.

\subsection{The $SU_q(2)$-coaction} 
We now give a coaction of the quantum group $SU_q(2)$ on 
the sphere
$S^7_q$. 
This coaction will be used later in Sect. \ref{se:hg} when analyzing the quantum principle bundle structure.
\\ Let us observe that the two pairs of generators $(\uno,\due),
(\tre,\qu)$ both yield a quantum plane,
\begin{eqnarray*}
\uno \due = q \due \uno \; ,  & \quad& \unob \dueb = q^{-1} \dueb \unob \; ,
\\
\tre \qu = q \qu \tre \; , && \treb \qub =q^{-1} \qub \treb \; .
\end{eqnarray*} Then we shall look for a right-coaction of $SU_q(2)$ 
on the rows of
the matrix $v$ in \eqref{u}. Other pairs of generators yield
quantum planes but
the only choice which gives a projection with the right number of 
generators is the one given above.

The  defining matrix of the quantum group $SU_q(2)$ reads
\begin{equation}\label{suq2}
\begin{pmatrix}
\alpha  & -q \bar{\gamma} \\
\gamma & \bar{\alpha}
\end{pmatrix}
\end{equation} with commutation relations
\cite{wo},
\begin{equation}\label{SU2}
\begin{array}{lll}
\alpha \gamma = q \gamma \alpha \; , \qquad & \alpha \bar{\gamma} = q
\bar{\gamma} \alpha \; , \qquad & \gamma 
\bar{\gamma}=\bar{\gamma}\gamma \; ,  \\
\alpha \bar{\alpha} + q^2 \bar{\gamma} \gamma = 1 \; , & \bar{\alpha} \alpha +
\bar{\gamma} \gamma = 1 \; .
\end{array}
\end{equation} We define a coaction of $SU_q(2)$ on the matrix \eqref{u} by,
\begin{equation}\label{cosu2}
\delta_R (v) :=
\left(
\begin{array}{cc} q^{-3} \uno & q^{-2} \due
\\- q^{-1}\dueb & q^{-1} \unob
\\
    q^{-1} \tre & -\qu
\\ - \qub & -\treb
\end{array}
\right)
\stackrel{.}{\otimes}
\begin{pmatrix}
\alpha  & -q \bar{\gamma} \\
\gamma & \bar{\alpha}
\end{pmatrix} \; .
\end{equation} We shall prove presently that this coaction comes from 
a coaction of
$A(SU_q(2))$ on the sphere algebra $A(S^7_q)$. For the moment we 
remark that, by
its form in
\eqref{cosu2} the entries of the projection $p=vv^*$ are automatically
coinvariants, a fact that we shall also prove explicitly in the following.

On the generators, the coaction \eqref{cosu2} is given explicitly by
\begin{equation} \label{coactionp}
\begin{array}{ll}
\delta_R(\uno)=  \uno \otimes \alpha + q ~\due \otimes \gamma \; ,  &
\delta_R(\unob) = q \dueb \otimes \bar{\gamma} + \unob \otimes
\bar{\alpha} =\overline{\delta_R(\uno)} \; ,
\\
\delta_R(\due)= - \uno \otimes \bar{\gamma} + \due \otimes \bar{\alpha} \; ,  &
\delta_R(\dueb) = \dueb \otimes \alpha - \unob \otimes \gamma
=\overline{\delta_R(\due)} \; ,
\\
\delta_R(\tre)=  \tre \otimes \alpha - q~\qu \otimes \gamma \; ,  &
\delta_R(\treb) =- q \qub \otimes \bar{\gamma} + \treb \otimes
\bar{\alpha}  =\overline{\delta_R(\tre)} \; ,
\\
\delta_R(\qu)=  \tre \otimes \bar{\gamma} + \qu \otimes \bar{\alpha} \; ,  &
\delta_R(\qub) = \qub \otimes \alpha + \treb \otimes  \gamma 
=\overline{\delta_R(\qu)}
\; ,
\end{array}
\end{equation} from which it is also clear its compatibility with the
anti-involution, i.e. $\delta_R(\bar{x}^i)=\overline{\delta_R(x_i)}$. The 
map $\delta_R$
in \eqref{coactionp} extends as an algebra homomorphism to the whole of
$A(S^{7}_q)$. Then, as alluded to before, we have the following
\begin{prop}\label{p4}
   The coaction \eqref{coactionp} is a right coaction of the quantum
group $SU_q(2)$
on the 7-sphere  $S^7_q$,
\begin{equation}
\delta_R : A(S^{7}_q) \lr A(S^{7}_q) \ot A(SU_q(2)) \; .
\end{equation}
\end{prop}
\noindent
\textit{Proof}.
    By using the commutation relations of $A(SU_q(2))$ in
\eqref{SU2}, a lengthy but easy computation gives that the 
commutation relations of
$A(S^7_q)$ are preserved. This fact also shows that extending $\delta_R$ as an
algebra homomorphism yields a consistent coaction.
\qed

\bigskip

\begin{prop}\label{equiv4} The algebra $A(S_q^4)$ is the algebra of 
coinvariants
under the coaction defined in \eqref{coactionp}.
\end{prop}
\noindent {\it Proof.} We have to show that $A(S^4_q)= \{ f \in A(S^7_q) ~|~
\delta_R(f)=f \ot 1 \}$.
    By using the commutation relations of $A(S_q^7)$ and those of
$A(SU_q(2))$, we
first prove explicitly that the generators of
$A(S_q^4)$ are coinvariants:
\begin{eqnarray*}
\delta_R(a)&=& q^{-4} \delta_R(\uno) \delta_R(\treb)-q^{-2} \delta_R(\due)
\delta_R(\qub)\\ &=& q^{-4}
\uno \treb \otimes
                (\alpha \bar{\alpha} + q^2 \bar{\gamma} \gamma)  -
q^{-2} \due \qub
\otimes (\gamma \bar{\gamma} + \bar{\alpha} \alpha)\\ &=& (q^{-4} \uno \treb -
q^{-2}\due \qub)
\otimes 1=  a \otimes 1
\\
\\
\delta_R(b)&=&- q^{-3}\delta_R (\uno) \delta_R(\qu) - q^{-2} \delta_R (\due)
\delta_R(\tre ) \\ &=& -q^{-3} \uno \qu \otimes  (\alpha \bar{\alpha} + q^2
\bar{\gamma}
\gamma) - q^{-2}
\due \tre
\otimes (\gamma \bar{\gamma} + \bar{\alpha} \alpha)
\\ &=&  - (q^{-3} \uno \qu  + q^{-2} \due \tre) \otimes 1  = b\otimes 1
\\
\\
\delta_R(t)&=&  q^{-2} \delta_R(\dueb) \delta_R(\due)+ q^{-2} \delta_R(\unob)
\delta_R(\uno)\\ &= & q^{-2} \dueb \due \otimes (\alpha \bar{\alpha} + q^2
\bar{\gamma}
\gamma)  + q^{-2} \unob \uno
\otimes (\gamma \bar{\gamma} + \bar{\alpha} \alpha)
\\ &=& (q^{-2} \dueb \due  + q^{-2} \unob \uno )\otimes 1 =  t\otimes 1
\end{eqnarray*} By construction the coaction is compatible  with the
anti-involution so that
$$
\delta_R(\bar{a})=\overline{\delta_R(a)}=\bar{a} \otimes 1 ,~~
\delta_R(\bar{b})=\overline{\delta_R(b)}=\bar{b} \otimes 1
$$ In fact, this only shows that $A(S_q^4)$ is made of coinvariants 
but does not
rule out the possibility of other coinvariants not in
$A(S_q^4)$. However this does not happen for the following reason.
From eq. \eqref{coactionp} it is clear that $w_1 \in \{\uno, \tre, \dueb, \qub \}$ (respectively $w_{-1} \in \{ \due, \qu, \unob, \treb \}$) are weight vectors of weight $1$
(resp. $-1$) in the fundamental comodule of $SU_q(2)$.
It follows that the only possible coinvariants are of the form $(w_1 w_{-1} -q w_{-1} w_1)^n$. When $n=1$ these are just the generators of
$A(S^4_q)$. \qed
\\
\begin{rema} The last part of the proof above is also related to the quantum Pl\"{u}cker coordinates. For every $2 \times 2$ matrix of \eqref{u}, let us define the determinant by
\begin{equation} det
\begin{pmatrix} a_{11} ~ & a_{12}
\\ a_{21} & a_{22}
\end{pmatrix} := a_{11} a_{22} - q ~ a_{12}a_{21} \; .
\end{equation} (Note that  $a_{12},~a_{21}$ do not commute and so in
the previous
formula the ordering between them is fixed.) Let $m_{ij}$ be the minors of
\eqref{u} obtained by considering the $i,j$ rows. Then
\begin{equation}\begin{array}{ll} m_{12}=q^{2}p_{11}= t \; ,   &
m_{13}=p_{14}= b
\; , \\ m_{14}=-q~p_{13} =-q ~a \; , & m_{23}=p_{24}=-q^2 \bar{a} \;
, \\ m_{24}=-q
~p_{23}= -q^{-1} \bar{b} \; , & m_{34}=-q~ p_{33}= q^{-3} t -q \; .
\end{array}
\end{equation} At $q=1$, these give the classical Pl\"{u}cker coordinates \cite{at}.
\end{rema}

\noindent The right coaction of $SU_q(2)$ on the 7-sphere  $S^7_q$ 
can be written as
\begin{equation}\label{coaction7}
\delta_R(\uno , \due ,\tre ,\qu)= (\uno , \due , \tre , \qu)
\stackrel{.}{\otimes}
\left(
\begin{array}{cccc}
\alpha ~~& -\bar{\gamma} ~~ &0 ~~ &0
\\ q \gamma & \bar{\alpha} & 0 &0
\\ 0 & 0 & \alpha &  \bar{\gamma}
\\ 0 &0 &-q \gamma &\bar{\alpha}
\end{array}
\right) \; ,
\end{equation} together with $\delta_R(\bar{x}_i)= 
\overline{\delta_R(x_i)}$. \\ In the
block-diagonal matrix which appears in  \eqref{coaction7} the second 
copy is given
by $SU_q(2)$ while the first one is twisted as
$$
\begin{pmatrix}
\alpha & -\bar{\gamma} \\ q \gamma &  \bar{\alpha}
\end{pmatrix}
    =
\begin{pmatrix} 1 & 0 \\ 0& - 1
\end{pmatrix}
\begin{pmatrix}
\alpha & \bar{\gamma} \\ -q \gamma &  \bar{\alpha}
\end{pmatrix}
\begin{pmatrix} 1 & 0 \\ 0& - 1
\end{pmatrix} \; .
$$ A similar phenomenon occurs in \cite{bct}.

\begin{rema} It is also interesting to observe that
$$
\delta_R(v^*v)=v^*v \otimes 1 = 1\otimes 1 \; .
$$ Indeed, \begin{eqnarray*}
\delta_R(\hs{\phi_1}{\phi_1})&=& \delta_R(q^{-6} \unob \uno +q^{-2} \due
\dueb + q^{-2} \treb \tre +
\qu \qub)
\\ &=& (-q^{-5} \dueb \uno + q^{-2} \uno \dueb + q^{-1} \qub \tre -
\tre \qub)\otimes
\bar{\gamma} \alpha
\\ && + (q^{-4} \dueb \due + q^{-2} \uno \unob + \qub \qu + \tre
\treb) \otimes \bar{\gamma}
\gamma
\\ && + (q^{-6} \unob \uno +q^{-2} \due \dueb + q^{-2} \treb \tre +
\qu \qub)\otimes
\bar{\alpha} \alpha
\\ && +
    (-q^{-5} \unob \due + q^{-2} \due \unob +q^{-1} \treb \qu - \qu
\treb)\otimes \bar{\alpha}
\gamma
\\ &=& \hs{\phi_2}{\phi_1} \otimes \bar{\gamma} \alpha +
\hs{\phi_2}{\phi_2} \otimes \bar{\gamma}
\gamma  +
\hs{\phi_1}{\phi_1} \otimes \bar{\alpha} \alpha + \hs{\phi_1}{\phi_2}
\otimes \bar{\alpha}
\gamma
\\ &=& 1 \otimes (\bar{\gamma} \gamma + \bar{\alpha} \alpha)= 1 \otimes 1 \; ,
\end{eqnarray*}
\begin{eqnarray*}
\delta_R(\hs{\phi_2}{\phi_2})&=& \delta_R (q^{-2} \uno \unob + q^{-4}
\dueb \due + \tre \treb + \qub
\qu)
\\ &=& (q^{-4} \dueb \uno -q^{-1}\uno \dueb - \qub \tre +  q \tre
\qub )\otimes \alpha
\bar{\gamma}
\\ && +
    (q^{-4} \dueb \due +q^{-2} \uno \unob + \qub \qu +  \tre \treb )
\otimes \alpha \bar{\alpha}
\\ && + (q^{-4} \unob \uno +  \due \dueb + \treb \tre  +q^2 \qu
\qub)\otimes \gamma \bar{\gamma}
\\ && + (q^{-4} \unob \due -q^{-1} \due \unob - \treb \qu   + q \qu
\treb) \otimes \gamma
\bar{\alpha}
\\ &=& -q \hs{\phi_2}{\phi_1} \otimes \alpha \bar{\gamma}+
\hs{\phi_2}{\phi_2} \otimes \alpha
\bar{\alpha} +
    q^2 \hs{\phi_1}{\phi_1} \otimes \gamma \bar{\gamma} -q
\hs{\phi_1}{\phi_2} \otimes \gamma
\bar{\alpha}
\\ &=& 1 \otimes (\alpha \bar{\alpha} + q^2 \gamma \bar{\gamma})= 1
\otimes 1 \; ,
\\
\\
\delta_R(\hs{\phi_1}{\phi_2})&=& q^{-5} \delta_R(\unob) \delta_R(\due) -q^{-2}
\delta_R(\due)
\delta_R(\unob)\nn
\\ 
&& \qquad \qquad- q^{-1}\delta_R(\treb)
\delta_R(\qu) + \delta_R(\qu) \delta_R(\treb) = 0
\end{eqnarray*} since $\delta_R$ defines a coaction on $S^7_q$ and so 
preserves its
commutation relations.\qed
\end{rema}

\section{Representations of the algebra $A(S^4_q)$}\label{se:rep} Let us now
construct irreducible $*$-representations  of
$A(S^4_q)$ as bounded operators on a separable Hilbert space $\hil$. For the
moment, we denote in the same way the elements of the algebra and 
their images as
operators in the given representation. As mentioned before, since 
$q\mapsto q^{-1}$
gives an isomorphic algebra, we can restrict ourselves to  $|q|<1$. We will
consider the representations which are $t$-finite \cite{KS}, i.e. such that the
eigenvectors of $t$ span $\hil$.

Since the self-adjoint operator
$t$ must be bounded due to the spherical relations, from the 
commutation relations
$ta= q^{-2} at,~~ t\bar{b}=q^{-4}\bar{b}t,$ it follows that the 
spectrum should be
of the form $\lambda q^{2k}$ and
$a, \bar b$ (resp. $\bar a, b$) act as rising (resp. lowering) operators on the
eigenvectors of $t$. Then boundedness implies the existence of a 
highest weight
vector, i.e. there exists a vector
$\ket{0,0}$ such that
\begin{equation} 
t\ket{0,0}=t_{00}\ket{0,0},~~~a\ket{0,0}=0,~~\bar{b}\ket{0,0}=0 ~.
\end{equation} By evaluating $q^4 \bar{a} a + b \bar{b}=(1-q^{-4}t)t$ 
on $\ket{0,0}$
we have
$$(1-q^{-4}t_{00})t_{00}=0$$ According to the values of the 
eigenvalue $t_{00}$ we
have two representations.

\subsection{The representation $\beta$}\label{se:pi} The first 
representation,
that we call $\beta,$ is obtained for
$t_{00}=0$. Then, $t\ket{0,0}=0$ implies
$t=0$. Moreover, using the commutation relations \eqref{s4} and \eqref{sr4}, it
follows that this representation  is the trivial one
\begin{equation} \label{pi} t=0,~~~a=0,~~~b=0 \; ,
\end{equation} the representation Hilbert space being just $\IC$; of 
course, $\beta
(1)=1$.

\subsection{The representation $\sigma$}\label{se:sigma} The second
representation, that we call $\sigma,$ is obtained for
$t_{00}=q^4$. This is infinite dimensional. 
We take the set
$\ket{m,n}= N_{mn} \bar{a}^m b^n \ket{0,0}$ with $n,m \in \IN$, to be an
orthonormal basis of the representation Hilbert space $\hil$, with 
$N_{00}=1$ and
$N_{mn}\in\IR$ the normalizations, to be computed below. \\ Then
$$
\begin{array}{l} t\ket{m,n}=t_{mn}\ket{m,n} ~,
\\
\bar{a}\ket{m,n}=a_{mn} \ket{m+1,n} ~,
\\ b\ket{m,n}=b_{mn} \ket{m,n+1} ~.
\end{array}
$$ By requiring that we have a $*$-representation we have also that
$$ a\ket{m,n}=a_{m-1,n}\ket{m-1,n}~, ~~ \bar{b}\ket{m,n}=b_{m,n-1}\ket{m,n-1} 
~,
$$ with the following recursion relations
$$ a_{m,n\pm 1} =q^{\pm 2} a_{m,n}~, \qquad b_{m \pm 1,n} = q^{\pm 2} b_{m,n}~,
\qquad b_{m,n} = q^2 a_{2n+1,m} ~.
$$ By explicit computation, we find
\begin{equation}
\begin{array}{l} t_{m,n}= q^{2m +4n+4} ~,
\\
\\ a_{m,n}=N_{mn} N_{m+1,n}^{-1} = (1-q^{2m+2})^{\frac{1}{2}}q^{m+2n+1}~,
\\
\\ b_{m,n}=N_{mn} N_{m,n+1}^{-1} =(1-q^{4n+4})^{\frac{1}{2}}q^{2(m+n+2)}~.
\end{array}
\end{equation}

\noindent In conclusion we have the following action
\begin{eqnarray} \label{sigma} && t\ket{m,n}= q^{2m +4n+4} \ket{m,n},
\\ && \bar a\ket{m,n}= (1-q^{2m+2})^{\frac{1}{2}}q^{m+2n+1}\ket{m+1,n},\nn \\
&&a\ket{m,n}= (1-q^{2m})^{\frac{1}{2}}q^{m+2n}\ket{m-1,n},\nn
\\ && 
b\ket{m,n}=(1-q^{4n+4})^{\frac{1}{2}}q^{2(m+n+2)}\ket{m,n+1},\nn \\ 
&& \bar
b\ket{m,n}=(1-q^{4n})^{\frac{1}{2}}q^{2(m+n+1)}\ket{m,n-1}.\nn
\end{eqnarray} It is straightforward to check that all the defining relations
\eqref{s4} and
\eqref{sr4} are satisfied.

\bigskip
\noindent In this representation the algebra generators are all trace class:
\begin{eqnarray}\label{classetraccia}
\Tr(t)&=& q^4 \sum_m q^{2m} \sum_n q^{4n}=\frac{q^4}{(1-q^2)(1-q^4)}
\; ~, \nn \\
\nn \\
\Tr(|a|)&=&q
\sum_{m,n}(1-q^{2m+2})^\frac{1}{2}q^{m+2n}=\frac{q}{1-q^2}\sum_{m}(1-q^{2m+2})^\frac{1}{2}q^{m}
\nn \\ &&\leq \frac{q}{1-q^2}  \sum_{m} q^{m}= \frac{q}{(1-q)(1-q^2)} \; ~, \\
\nn \\
\Tr(|b|)&=&q^4
\sum_{m,n}(1-q^{4n+4})^\frac{1}{2}q^{2(n+m)} =
\frac{q^4}{1-q^2}\sum_{n}(1-q^{4n+4})^\frac{1}{2}q^{2n}
\nn \\ &&\leq \frac{q^4}{1-q^2}  \sum_{n} q^{2n}= 
\frac{q^4}{(1-q^2)^2} \; ~. \nn
\end{eqnarray}

 From the sequence of Schatten ideals in the algebra of compact 
operators one knows
\cite{si} that the norm closure of trace class operators gives the 
ideal of compact
operators $\ck$. As a consequence, the closure of $A(S_q^4)$ is the 
$C^*$-algebra
$\cc(S_q^4)=\ck \op \IC \II$.

\section{The index pairings}\label{in-pa}

The `defining' self-adjoint idempotent $p$ in \eqref{proj} determines 
a class in
the $K$-theory of $S_q^4$, i.e. $[p]\in K_0[\cc(S_q^4)]$. A way to prove its
nontriviality is by pairing it with a nontrivial element in the dual 
$K$-homology,
i.e. with (the class of) a nontrivial Fredholm module $[\mu]\in 
K^0[\cc(S_q^4)]$.
In fact, in order to compute the pairing of
$K$-theory with $K$-homology, it is more convenient to first compute the
corresponding Chern characters in the cyclic homology
$\ch_*(p) \in HC_*[A(S^{4}_q)]$ and cyclic cohomology $\ch^*(\mu)\in
HC^*[A(S^{4}_q)]$ respectively, and then use the pairing between 
cyclic homology
and cohomology \cite{co}.

Like it happens for the $q$-monopole \cite{HM}, to compute the pairing and to prove the nontriviality of the
bundle it is enough to consider $HC_0[A(S^{4}_q)]$ and dually to take a suitable trace of the projector. 

\bigskip 
The Chern character of the projection $p$ in \eqref{proj} has a component in degree zero
$\ch_0(p)\in HC_0[A(S^{4}_q)]$ simply given by the matrix trace,
\begin{equation}\label{chp}
\ch_0(p) := \tr(p) = 2 - q^{-4} (1-q^2)(1-q^4) ~t ~\in A(S^4_q).
\end{equation} The higher degree parts of $\ch_*(p)$ are obtained via the
periodicity operator $S$; not needing them here we shall not dwell 
more upon this
point and refer to
\cite{co} for the relevant details.

\bigskip As mentioned, the K-homology of an involutive algebra $\ca$ 
is given in
terms of homotopy classes of Fredholm modules. In the present situation we are
dealing with a
$1$-summable Fredholm module $[\mu]\in K^0[\cc(S_q^4)]$. This is in 
contrast to the
fact that the analogous element of $K_0(S^4)$ for the undeformed 
sphere is given by
a
$4$-summable Fredholm module, being the fundamental class of  $S^4$.

The Fredholm module $\mu := (\hil,\Psi,\gamma)$ is constructed as follows. The 
Hilbert space
is $\hil=\hil_\sigma \oplus \hil_\sigma$ and the representation is
$\Psi=\sigma \oplus \beta$. Here $\sigma$ is the representation of
$A(S^4_q)$ introduced in
\eqref{sigma} and $\beta$ given in \eqref{pi} is trivially extended to
$\hil_\sigma$. The grading operator is
\[
\gamma = \begin{pmatrix} 1&0
\\0&-1
\end{pmatrix} .
\] The corresponding Chern character $\ch^*(\mu)$ of the class of this Fredholm
module has  a component in degree $0$, $\ch^{0}(\mu) \in HC^0[A(S^{2n}_q)]$.
   From the general construction \cite{co}, the element
$\ch^0(\mu_{\mathrm{ev}})$ is the trace
\begin{equation}\label{tau1}
\tau^1(x) := \Tr\left(\gamma \Psi(x)\right) =
\Tr\left(\sigma(x) - \beta(x)\right).
\end{equation} The operator $\sigma(x) - \beta(x)$ is always trace 
class. Obviously
$\tau^1(1)=0$. The higher degree parts of $\ch^*(\mu_{\mathrm{ev}})$ 
can again be
obtained via a periodicity operator.
\\
A similar construction of the class $[\mu]$ and the corresponding Chern character were given in \cite{MNW} for quantum two and three 
dimensional spheres.

\bigskip We are ready to compute the pairing:
\begin{align}
\left\langle[\mu],[p]\right\rangle &:= \left\langle \ch^0(\mu),
\ch_0(p)\right\rangle = -q^{-4}(1-q^2)(1-q^4) ~\tau^1(t) \nn
\\  &= -q^{-4}(1-q^2)(1-q^4) \Tr(t) = -q^{-4}(1-q^2)(1-q^4)
q^{4}(1-q^2)^{-1}(1-q^4)^{-1} \nn \\ &= - 1 ~.
\label{pair}
\end{align} This result shows also that the right $A(S^4_q)$-module
$p[A(S_q^4)^4]$ is not free. Indeed, any free module is represented in
$K_0[\cc(S_q^4)]$ by the idempotent $1$, and since
$\left\langle[\mu],[1]\right\rangle=0$, the evaluation of $[\mu]$ on any free
module always gives zero.

\bigskip We can extract the `trivial' element in the $K$-homology
$K^0[\cc(S_q^4)]$ of the quantum sphere $S_q^4$ and use it to measure 
the `rank' of
the idempotent
$p$. This generator corresponds to the trivial generator of the $K$-homology
$K_0(S^4)$ of the classical sphere
$S^4$. The latter (classical) generator is the image of the generator of the
$K$-homology of a point by the functorial map $K_*(\iota) : K_0(*) \to
K_0(S^{N})$, where $\iota : *
\hookrightarrow S^{N}$ is the inclusion of a point into the sphere. Now, the
quantum sphere
$S_q^4$ has just one `classical point', i.e. the $1$-dimensional representation
$\beta$ constructed in Sect.~\ref{se:pi}. The corresponding
$1$-summable Fredholm module $[\varepsilon]\in K^0[\cc(S_q^4)]$   is easily
described: the Hilbert space is $\IC$ with representation $\beta$; the grading
operator is
$\gamma=1$. Then the degree $0$ component $\ch^{0}(\varepsilon) \in
HC^0[A(S^{2n}_q)]$ of the corresponding Chern character is the trace 
given by the
representation itself (since it is a homomorphism to a commutative algebra),
\begin{equation}
\tau^0(x) = \beta (x) ~,
\end{equation} and vanishes on all the generators whereas $\tau^0(1)=1$. \\ Not
surprisingly, the pairing with the class of the idempotent $p$ is
\begin{equation}
\left\langle[\varepsilon],[p]\right\rangle := \tau^0(\ch_0(p)) =
\beta (2) = 2 ~.
\end{equation}

\section{Quantum principal bundle structure}\label{se:hg} 

Recall that 
if $H$ is a Hopf algebra and $P$ a right $H$-comodule 
algebra with
multiplication $m: P\ot P \raw P$ and  coaction $\Delta_R : P \raw P 
\ot H$ and $B\subseteq P$ is the subalgebra of coinvariants, the extension $B\subseteq 
P$ is $H$ Hopf-Galois if the canonical map
\begin{equation}\label{canbis}
\chi : P \otimes_B P \lr P \ot H \; , 
p' \ot_B p
\mapsto
\chi(p' \otimes_B p) = p' p_{(0)} \ot p_{(1)} \; ,
\end{equation} is bijective.
As mentioned, for us a quantum principle bundle will be the same as a Hopf-Galois extension.
For quantum structure groups which are cosemisimple and have bijective antipodes, as is the case for $SU_q(2)$,
Th.~I of \cite{scH}
grants further nice properties. In particular the surjectivity of the canonical map implies bijectivity and faithfully flatness
of the extension.
Moreover, an additional useful result \cite{scP} is that the map $\chi$ is
surjective whenever, for any generator $h$ of $H$, the element $1
\ot h$ is in its image. This follows from the  left $P$-linearity and right
$H$-colinearity of the map $\chi$. Indeed, let $h,~k$ be two elements
of $H$ and
$\sum p_i' \ot p_i,~\sum q_j' \ot q_j \in P \ot P$ be such that
$\chi(\sum p_i' \ot_B p_i)=1 \ot h,~\chi(\sum q_j' \ot_B q_j)=1 \ot k$. Then
$\chi(\sum p_i'q_j' \ot_B q_j p_i)=1 \ot kh$, that is $1 \ot kh$ is 
in the image of
$\chi$. But, since the map $\chi$ is left  $P$-linear, this implies its
surjectivity.

\bigskip

\begin{defi}\label{de}  Let $P$ be a bimodule over the ring $B$. Given any two
elements $\ket{\xi_1}$ and $\ket{\xi_2}$ in the free module
$\ce=\IC^m\ot P$, we shall define $\hstb{\xi_1}{\xi_2}
\in P \otimes_B P$ by
\begin{equation}
\hstb{\xi_1}{\xi_2}:= \sum_{j=1}^m \bar\xi_1^j \ot_B \xi_2^j \;.
\end{equation} Analogously, one can define quantities 
$\hst{\xi_1}{\xi_2} \in P \ot
P$ with the same formula as above and tensor products taken over the 
ground field
$\IC$.
\end{defi}

\begin{prop} The extension $A(S^{7}_q) \subset  A(Sp_q(2))$ is a 
faithfully flat
$A(Sp_q(1))$-Hopf-Galois extension.
\end{prop}

\noindent {\it Proof.} Now $P=A(Sp_q(2))$, $H=A(Sp_q(1))$ and
$B=A(S^{7}_q)$ and the
coaction
$\Delta_R$ of
$H$ is given just before Prop.~\ref{p2}. Since $A(Sp_q(1))\simeq 
A(SU_{q^2}(2))$
has a bijective antipode and is cosemisimple (\cite{KS}, Chapter 11), from the
general considerations given above in order to show the bijectivity of the
canonical map
$$
\chi : A(Sp_q(2)) \otimes_{ A(S^{7}_q) } A(Sp_q(2))
\lr  A(Sp_q(2)) \ot A(Sp_q(1)) \; ,
$$ it is enough to show that all  generators $\alpha, \gamma, \bar\alpha,
\bar\gamma $ of $A(Sp_q(1))$  in \eqref{spq1} are in its image. \\ Let
$\ket{T^2},\ket{T^3}$ be the  second and third columns of the 
defining matrix $T$
of $Sp_q(2)$. We  shall think of them as elements  of the free module 
$\IC ^4\ot
A(Sp_q(2))$. Obviously,
$\hs{T^i}{T^j}=\delta^{ij}$.  Recalling that
$A(Sp_q(2))$ is both a left and right $A(S^7_q)$-module and using 
Def.~\ref{de},
we have  that
$$
\chi
\left(\begin{array}{ll}
\hsts{T^2}{T^2} & \hsts{T^2}{T^3}
\\ \\
\hsts{T^3}{T^2} & \hsts{T^3}{T^3}
\end{array}
\right) = 1
\pot
\left(
\begin{array}{lr}
\alpha & -q^2 \bar\gamma \\ \\
\gamma &
\bar\alpha
\end{array}
\right)\; .$$ Indeed,
\begin{eqnarray*}
\chi  (\hsts{T^2}{T^2}) &=& \overline T^2_i \Delta_R T^2_i=
\hs{T^2}{T^2}\ot
\alpha+\hs{T^2}{T^3}\ot \gamma=1\ot \alpha  \; , \\
\chi (\hsts{T^3}{T^2})  &=& \overline T^3_i \Delta_R T^2_i=
\hs{T^3}{T^2}\ot
\alpha+\hs{T^3}{T^3}\ot \gamma=1\ot \gamma ;
\end{eqnarray*} a similar computation giving the other two  generators.
\qed

\begin{prop}\label{p7} The extension $A(S^{4}_q)
\subset  A(S^{7}_q)$ is a faithfully flat
$A(SU_q(2))$-Hopf-Galois  extension.
\end{prop}
\noindent {\it Proof.} Now $P=A(S^{7}_q)$,
$H=A(SU_q(2))$ and $B=A(S^{4}_q)$ and the coaction $\delta_R$ of $H$
is given  in
Prop.~\ref{p4}. As already mentioned $A(SU_q(2))$ has a bijective 
antipode and is
cosemisimple, then as before  in order to show the  bijectivity of 
the canonical map
$$
\chi : A(S^{7}_q) \otimes_{  A(S^{4}_q) } A(S^{7}_q)
\lr  A(S^{7}_q) \ot A(SU_q(2)) \; ,
$$ we  have to show that all generators $\alpha, \gamma,
\bar\alpha,
\bar\gamma $ of $A(SU_q(2))$ in \eqref{suq2} are in its  image. \\ 
Recalling that
$A(S^7_q)$ is both a left and right
$A(S^4_q)$-module and using Def.~\ref{de}, we have  that
$$
\chi
\left(\begin{array}{ll}
\hstq{\phi_1}{\phi_1} &
\hstq{\phi_1}{\phi_2} \\ \\
\hstq{\phi_2}{\phi_1} &
\hstq{\phi_2}{\phi_2}
\end{array}
\right) = 1
\pot
\left(
\begin{array}{lr}
\alpha & -q \bar\gamma \\ \\
\gamma &
\bar\alpha
\end{array}
\right)\; ,$$ where
$\ket{\phi_1},\ket{\phi_2}$ are the two vectors introduced  in 
eqs.~\eqref{psi1}
and \eqref{psi2}.  Indeed

\begin{eqnarray*}
\chi(\hstq{\phi_1}{\phi_1}) &=&
\chi\left(q^{-6} \unob \ota \uno + q^{-2}
\due \ota \dueb \right.
\\ && \qquad \qquad \left. + q^{-2} \treb \ota \tre + \qu \ota
\qub
\right) \\ &=&  q^{-6} \unob \delta_R( \uno) + q^{-2} \due \delta_R(
\dueb) + q^{-2} \treb
\delta_R(\tre) + \qu \delta_R( \qub) \\ &=& q^{-6}
\unob \uno \ot \alpha +q^{-5} \unob \due \ot \gamma + q^{-2}
\due
\dueb \ot \alpha -q^{-2} \due \unob \ot \gamma \\ && + q^{-2}
\treb \tre \ot \alpha -q^{-1} \treb \qu \ot \gamma + \qu \qub
\ot
\alpha + \qu \treb \ot \gamma \\ &=& \hs{\phi_1}{\phi_1} \ot \alpha = 
1 \ot \alpha
\;  ,
\end{eqnarray*}
\begin{eqnarray*}
\chi(\hstq{\phi_2}{\phi_1}) &=&  q^{-5} \dueb \delta_R (\uno) -q^{-2} \uno
\delta_R (\dueb) -q^{-1} \qub
\delta_R (\tre) + \tre \delta_R (\qub) \\ &=& q^{-5} \dueb \uno \ot
\alpha + q^{-4} \dueb \due \ot \gamma -q^{-2} \uno
\dueb \ot \alpha  +q^{-2} \uno \unob \ot \gamma \\ && - q^{-1} \qub 
\tre \ot \alpha
+ \qub \qu \ot \gamma + \tre \qub \ot
\alpha + \tre \treb \ot \gamma \\ &=& \hs{\phi_2}{\phi_1} \ot \alpha +
\hs{\phi_2}{\phi_2} \ot \gamma = 1 \ot
\gamma \; ,
\end{eqnarray*} with similar computations for the other generators. \qed

It was proven in \cite{bcdt} that the bundle constructed in
\cite{bct} is a coalgebra Galois extension \cite{BM98,BH}. The fact 
that our bundle
$A(S^4_q)\subset A(S^7_q)$  is Hopf-Galois shows also that these two 
bundles cannot be the same.

\bigskip On our extension $A(S^4_q) \subset A(S^7_q)$ there is a 
strong connection.
Indeed a $H$-Hopf-Galois extension $B \subseteq P$ for which $H$ is
cosemisimple
and has a bijective antipode is also  equivariantly projective, that is there
exists a left
$B$-linear right $H$-colinear splitting $s:P \raw B \ot P$ of the 
multiplication map
$m: B \ot P \raw P,~~m\circ s=id_P$ \cite{scsc}. Such a map 
characterizes the so
called strong connection. Constructing a strong connection is an alternative way to
prove that one has a Hopf Galois extension \cite{dgh,hP}.

In particular, if $H$ has an invertible antipode $S$, an
equivalent  description
of a strong connection can be given in terms of a map
$\ell: H \raw P \ot P$ satisfying a list of conditions \cite{maj, BH04} (see also
\cite{HMS,BDZ}). We denote by $\Delta$ the coproduct on $H$ with 
Sweedler notation
$\Delta(h)= h_{(1)}\ot h_{(2)}$, by $\delta:P\raw P \ot H$ the right-comodule
structure on
$P$ with notation $\delta p=p_{(0)} \ot p_{(1)}$, and
$\delta_l:P \raw H \ot P$ is the induced left $H$-comodule structure 
of P defined
by
$\delta_l (p) = S^{-1} (p_{(1)}) \ot p_{(0)}$. Then, for the map $\ell$ one
requires that
$\ell(1) =1 \ot 1$  and that  for all $h \in H$,
\begin{eqnarray}\label{ell} && \chi(\ell(h))=1 \ot h \, ,\nn \\ 
&&\ell(h_{(1)})\ot
h_{(2)} = (id \ot \delta)\circ \ell(h) \, ,\nn \\ &&h_{(1)}\ot 
\ell(h_{(2)}) = (
\delta_l \ot id)\circ \ell \, (h)
\end{eqnarray} The splitting $s$ of the multiplication map is then given by
\[ s: P \rightarrow  B\ot P \, ,\quad p \mapsto p_{(0)} \ell (p_{(1)})~.
\] Now, if $g,h \in H$ are such that $\ell(g)=g^{1}\ot g^{2}$ and $\ell(h)=
h^{1}\ot h^{2}$ satisfy condition \eqref{ell} so does
$\ell(gh)$ defined by
\begin{equation}\label{proell}
\ell(gh):=h^{1}g^{1}\ot g^{2}h^{2} \, .
\end{equation} If $H$ has a PBW basis \cite{ka},  this fact can be used to
iteratively construct $\ell$ once one knows its  value on the 
generators of $H$.

For $H=A(SU_q(2))$, with generators, $\alpha, \gamma, \bar\alpha$ and
$\bar\gamma$,
the PBW basis is given by $\alpha^k \gamma^l \bar\gamma^m$,  with $k,l,m\in
\{0,1,2, \dots \}$ and
$\gamma^k \bar\gamma^l \bar\alpha^m$,  with $k,l \in \{0,1,2, \dots 
\}$ and $m\in
\{1,2, \dots \}$ \cite{wo}. Then, for our extension $A(S^4_q) \subset 
A(S^7_q)$ the
map $\ell$ can be  constructed as follows. Firstly, we put $\ell(1) =1 \ot 1$.
Then, on the generators we set
\begin{eqnarray*} &\ell (\alpha) := \hst{\phi_1}{\phi_1} \, , \quad &\ell (\bar
\alpha) := \hst{\phi_2}{\phi_2} \, , \\ &\ell (\gamma) := 
\hst{\phi_2}{\phi_1} \, ,
\quad & \ell (\bar
\gamma) := -q^{-1}\hst{\phi_1}{\phi_2}
\, .
\end{eqnarray*} These expressions for $\ell$ satisfy all the 
properties \eqref{ell}:

\noindent Firstly, $\chi(\ell(\alpha))=1 \ot \alpha$ follows from the proof of
Prop.~\ref{p7}. Then,
\begin{eqnarray*} (id \ot \delta)\circ ~\ell(\alpha)&=& q^{-6}\unob 
\ot \delta \uno
+  q^{-2}\due \ot \delta \dueb + q^{-2}\treb \ot \delta \tre + \qu 
\ot \delta \qub
\\ &=&
\hst{\phi_1}{\phi_1} \ot \alpha + \hst{\phi_1}{\phi_2} \ot \gamma \\ &=&
\ell(\alpha) \ot \alpha - q \ell( \bar \gamma) \ot \gamma =
\ell(\alpha_{(1)})\ot \alpha_{(2)} \;  .
\end{eqnarray*} Moreover
\begin{eqnarray*} (\delta_l \ot id)\circ  ~\ell(\alpha)&=& q^{-6} 
(\alpha \ot \unob
-q^2 \bar \gamma \ot \dueb)
\ot \uno +  q^{-2} (q \bar \gamma \ot \uno + \alpha \ot \due) \ot
\dueb
\\ &+& q^{-2} (q^2 \bar \gamma \ot \qub + \alpha \ot \treb) \ot  \tre 
+ (-q \bar
\gamma \ot \tre + \alpha \ot \qu) \ot \qub
\\ &=&
\alpha \ot \hst{\phi_1}{\phi_1} - q \bar \gamma \ot 
\hst{\phi_2}{\phi_1} \\ &=&
\alpha \ot \ell(\alpha) - q  \bar \gamma \ot \ell(\gamma) =
\alpha_{(1)} \ot \ell(\alpha_{(2)}) \; .
\end{eqnarray*} Similar computations can be carried for $\gamma, 
\bar\alpha$ and
$\bar \gamma$.

That an iterative procedure constructed by using \eqref{proell} on 
the PBW  basis
leads to a well defined $\ell$ on the whole of  $H=A(SU_q(2))$ will 
be proven  in
the forthcoming paper \cite{future} where other elaborations coming from  the
existence of a strong connection will be presented as well.

\subsection{The associated bundle and the coequivariant maps}

We now give some elements of the theory of associated quantum vector bundles \cite{BM93} (see also
\cite{du}). Let $B\subset P$ be a $H$-Galois extension with $\Delta_R$ the
coaction of $H$ on $P$. Let $\rho : V \raw H \ot V$ be a
corepresentation of $H$
with $V$ a finite dimensional vector space.  A coequivariant map is an element $\varphi$ in $P \ot V$
with the property that
\begin{equation}
(\Delta_R \ot id)\varphi = (id \ot (S \ot id) \circ \rho)) \varphi  \, ;
\end{equation} where $S$ is the antipode of $H$.  The collection
$\Gamma_\rho(P,V)$ of coequivariant maps is a right and left $B$-module.

\bigskip The algebraic analogue of bundle nontriviality is translated 
in the fact
that the Hopf-Galois extension $B\subset P$ is not cleft.
On the  other
hand, it is  knows that for a cleft Hopf-Galois extension, the module of
coequivariant maps $\Gamma_\rho(P, V)$ is isomorphic to the
free module of coinvariant maps $\Gamma_0(P,V)= B \ot V$ \cite{BM93,HM}.

\bigskip

For our $A(SU_q(2))$-Hopf-Galois extension
$A(S^4_q)\subset A(S^7_q)$,  let $\rho_1 : \IC^2 \raw \IC^2 \ot
A(SU_q(2))$ be the
fundamental  corepresentation of $A(SU_q(2))$  with 
$\Gamma_1(A(S^7_q),\IC^2)$
the right $A(S^4_q)$-module of corresponding coequivariant maps.

Now, the projection $p$ in \eqref{proj} determines a quantum vector bundle over
$S_q^4$ whose  module of section is $p[A(S_q^4)^4]$, which is clearly a right
$A(S^4_q)$-module. The following proposition in straightforward
\begin{prop}  The modules ${\cal E}:=p [A(S_q^4)^4]$ and $\Gamma_1(A(S^7_q),\IC^2)$  are
isomorphic as right
$A(S^4_q)$-modules.
\end{prop}
\noindent
{\it Proof.} Remember that $p=vv^*$
with $v$ in \eqref{u}.
The element $p(F) \in {\cal E}$, with $F=(f_1,f_2,f_3,f_4)^t$,  corresponds to the equivariant map $v^*F \in \Gamma_1(A(S^7_q),\IC^2)$.
\qed
\\

We expect that a similar construction extends to every irreducible
corepresentation of
$A(SU_q(2))$ by means of suitable projections giving the
corresponding associated bundles \cite{future}.

\begin{prop} The Hopf-Galois extension
$A(S^4_q)\subset A(S^7_q)$ is not cleft.
\end{prop}

\noindent {\it Proof.} As mentioned, the cleftness of the extension 
does imply that
all modules of coequivariant maps are free. On the other hand,  the 
nontriviality
of the pairing \eqref{pair} between the defining  projection $p$ in 
\eqref{proj}
and the Fredholm module $\mu$ constructed  in Sect.~\ref{in-pa} also 
shows that the
module
$p [A(S_q^4)^4] \simeq \Gamma_\rho(A(S^7_q),\IC^2)$ is not free. \qed
\\

\subsection*{Acknowledgments} We thank the referee for many useful comments and suggestions. We are grateful to Tomasz
Brzezi{\'n}ski and Piotr M.
Hajac for several important remarks on a previous version of the 
compuscript. Also,
Eli Hawkins, Walter van Suijlekom, Marco Tarlini are thanked for very useful
discussions.

\appendix

\section{The classical Hopf fibration $S^7 \rightarrow S^4$}
\label{se:chf} We shall review the classical construction of the  basic
anti-instanton  bundle over the four dimensional sphere $S^4$  in a 
`noncommutative
parlance' following \cite{la00}. This has been  useful in the main text for our
construction of the quantum  deformation of the Hopf bundle.

We write the generic element of the  group $SU(2)$ as
\begin{equation} w =
\left(
\begin{array}{cc} w_1 &  w_2 \\ -\bar{w}_2 &
\bar{w}_1
\end{array}
\right)~.
\end{equation} The $SU(2)$ principal  fibration $SU(2) \raw S^7 \raw 
S^4$ over the
sphere $S^4$ is  explicitly realized as follows. The total space  is
$S^7 = \{z=(z_1, z_2, z_3, z_4) \in \IC^4 ~,~
\sum_{i=1}^4|z_i|^2 = 1 \}$~, with right diagonal action
\begin{equation}\label{sp1act} S^7 \times SU(2) ~\raw~  S^7~, ~~~z
\cdot w := (z_1,
z_2, z_3, z_4)
\left(
\begin{array}{cccc} w_1 ~~& w_2  ~~ &0 ~~ &0
\\ -\bar{w}_2 &
\bar{w}_1 & 0 &0
\\ 0 & 0 & w_1 ~~& w_2
\\ 0 &0 &-\bar{w}_2 &
\bar{w}_1
\end{array}
\right) \; .
\end{equation} The bundle  projection $\pi : S^7 \raw S^4$ is just the Hopf
projection and  it can be explicitly given as $\pi(z_1,z_2,z_3,z_4) :=
(x,\alpha,\beta)$ with
\begin{eqnarray}\label{s4coord} && x = |z_1|^2  + |z_2|^2 - |z_3|^2 - 
|z_4|^2 = - 1
+ 2 (|z_1|^2 +  |z_2|^2) = 1 -  2(|z_3|^2 + |z_4|^2) ~, \nn \\ && 
\alpha =2( z_1
\bar{z}_3 + z_2
\bar{z}_4) ~,  ~~~~~~\beta =2(- z_1  z_4 + z_2 z_3)  ~.
\end{eqnarray} One checks that $ |\alpha|^2 + |\beta|^2 + x^2 =
(\sum_{i=1}^4|z_i|^2)^2 = 1$~.

\bigskip

We need the rank $2$  complex vector bundle $E$ associated with the 
defining  left
representation $\rho$ of $SU(2)$ on $\IC^2$. The quickest way to  get 
this is to
identify $S^7$ with the unit sphere in the
$2$-dimensional  quaternionic (right) $\IH$-module
$\IH ^2$ and $S^4$  with the projective line $\IP^1(\IH)$, i.e. the set
   of equivalence  classes $(w_1,w_2)^t\simeq (w_1,w_2)^t \lambda$ with
$(w_1,w_2)\in  S^7$ and $\lambda \in Sp(1)\simeq SU(2)$. Identifying
$\mathbb{H}\simeq\IC ^2$, the vector
$(w_1,w_2)^t \in  S^7$ reads
\begin{equation}\label{ati} v=\left(
\begin{array}{cc} z_1&z_2\\ -\bar{z}_2&\bar{z}_1\\ z_3&z_4\\ 
-\bar{z}_4&\bar{z}_3
\end{array}
\right).
\end{equation} This  is actually a map from $S^7$ to the Stieffel 
variety of frames
for
$E$. In particular, notice that the two vectors
$\ket{\psi_1},\ket{\psi_2}$ given by the columns  of
$v$ are  orthonormal, indeed $v^*v=\II_2$. As a  consequence,
$ p :=vv^*= \ket{\psi_1}
\bra{\psi_1} + \ket{\psi_2} \bra{\psi_2}$ is a  self-adjoint idempotent (a projector), $p^2 = 
p$, $p^* = p$.
Of  course $p$ is $SU(2)$ invariant and hence its entries are 
functions  on $S^4$
rather than $S^7$. An explicit  computation yields
\begin{equation}\label{class} p =
\frac{1}{2}
\left(
\begin{array}{cccc} 1 + x & 0       & \alpha  & \beta \\ 0       &1 + x & -
\bar{\beta} & \bar{\alpha}
\\
\bar{\alpha} & -\beta &1 - x & 0 \\
\bar{\beta} & \alpha & 0 & 1 -  x
\end{array}
\right)~,
\end{equation} where $(x,\alpha,\beta)$ are  the coordinates (\ref{s4coord}) on
$S^4$. Then
$p \in
\Mat_4(C^\infty(S^4, \IC))$  is of rank $2$ by  construction. \\

The matrix $v$ in \eqref{ati} is a  particular example of the 
matrices $v$ given
in  \cite{at}, for
$n=1,~k=1,~C_0=0,~C_1=1,~ D_0=1,~D_1=0$. This gives the 
(anti-)instanton of charge
$-1$ centered at the origin and with unit  scale. The only difference 
is  that here
we identify $\IC ^4$ with
$\mathbb{H}^2$ as a right $\mathbb{H}$-module. This notwithstanding,  the
projections constructed in the two formalisms actually  coincide. 
Finally recall
that, as mentioned already, the  classical limit of our quantum projection
\eqref{proj} is conjugate  to \eqref{class}.

\bigskip The canonical connection  associated with the  projector,
\begin{equation}
\nabla := p \circ d ~:~
\Gamma^\infty(S^4,E) ~\raw~
\Gamma^\infty(S^4,E)
\otimes_{C^\infty(S^4, \IC)} \Omega^1(S^4,
\IC),
\end{equation} corresponds to a Lie-algebra valued ($su(2)$)
$1$-form $A$ on $S^7$  whose matrix components are given  by
\begin{equation}\label{connection} 
A_{ij} = \hs{\psi_i}{d \psi_j} \; , \qquad  i,j=1,2 \; .
\end{equation} 
This connection can be used to compute  the Chern 
character of the
bundle. Out of the curvature of the  connection
$\nabla^2 = p (dp)^2$ one has the Chern $2$-form and
$4$-form given  respectively  by
\begin{eqnarray}\label{inscf} &&C_1(p) := - \frac{1}{2 \pi i} ~\tr 
(p (dp)^2)~,
\nn \\ &&C_2(p) := - \frac{1}{8 \pi^2 } ~[ \tr (p  (dp)^4) - C_1(p)C_1(p) ]~,
\end{eqnarray} with the trace $\tr$ just  an ordinary matrix trace. 
It turns out
that the $2$-form $p (dp)^2$  has vanishing  trace so that $ C_1(p) =  0$. As for the second Chern class, a straightforward 
calculation shows
that,
\begin{eqnarray}\label{inscf2} C_2(p) &=&  -\frac{1}{32\pi^2} [(x_0 dx_4 - x_4
dx_0) (d\xi)^3 + 3 dx_0 dx_4 ~\xi  ~(d\xi)^2] \nn \\ &=& -\frac{3}{8\pi^2} [x_0
dx_1 dx_2 dx_3 dx_4 +  \mbox{\it cyclic permutations}]\nn \\ &=&  -\frac{3}{8\pi^2} ~d (vol(S^4))~.
\end{eqnarray} The second Chern  number is then given by
\begin{equation}\label{inscn2} c_2(p) =
\int_{S^4} C_2(p) = -\frac{3}{8\pi^2} \int_{S^4} d (vol(S^4))  =
-\frac{3}{8\pi^2}
\frac{8}{3} \pi^2 = -1~.
\end{equation} The  connection $A$ in \eqref{connection} is 
(anti-)self-dual, i.e.
its curvature
$F_A := d A + A
\wedge A$ satisfies (anti-)self-duality equations, $*_H F_A  = - 
F_A$, with $*_H$
the Hodge map of  the canonical (round) metric on  the sphere $S^4$. 
It is indeed
the basic Yang-Mills  anti-instanton  found in
\cite{bpst}.


\addcontentsline{toc}{section}{References}
\begin{thebibliography}{9}

\bibitem{at} M. Atiyah,
\emph{The geometry of Yang-Mills fields}, Lezioni Fermiane. Accademia
Nazionale dei
Lincei e Scuola Normale  Superiore, Pisa 1979.

\bibitem{bpst} A. Belavin, A. Polyakov, A. Schwartz, Y. Tyupkin, {\it
Pseudoparticles solutions of the Yang-Mills equations}, Phys. Lett.
58 B (1975)
85-87.


\bibitem{bct} F. Bonechi, N. Ciccoli, M. Tarlini,
\emph{Noncommutative instantons on the 4-sphere from quantum groups},
Commun. Math.
Phys. 226 (2002) 419-432.

\bibitem{bcdt} F. Bonechi, N. Ciccoli, L. D\c{a}browski, L. M. Tarlini,
\emph{Bijectivity of the canonical map for the non-commutative
instanton bundle},
J. Geom. Phys. 51 (2004) 71-81.

\bibitem{BDZ} T. Brzezi{\'n}ski, L. D\c{a}browski, B. Zielinski,
\emph{Hopf fibration and monopole connection over the contact quantum 
spheres}, J.
Geom. Phys. 50 (2004) 345-359.

\bibitem{BH} T. Brzezi{\'n}ski, P. M. Hajac,
\emph{Coalgebra extensions and algebra coextensions of Galois type}, Commun.
Algebra 27 (1999) 1347-1368.

\bibitem{BH04} T. Brzezi{\'n}ski, P. M. Hajac,
\emph{The Chern-Galois character}, C. R. Acad. Sci. Paris, Ser. I 333 (2004)
113-116.

\bibitem{BM93} T. Brzezi{\'n}ski, S. Majid,
\emph{Quantum group gauge theory on quantum spaces}, Commun. Math. 
Phys. 157 (1993)
591-638. Erratum 167 (1995) 235.

\bibitem{BM98} T. Brzezi{\'n}ski, S. Majid,
\emph{Coalgebra Bundles}, Commun. Math. Phys. 191 (1998) 467-492.

\bibitem{co} A. Connes,
\emph{Noncommutative geometry}, Academic Press, 1994.

\bibitem{du} M. Durdevich,  {\it Geometry of quantum principal 
bundles I}, Commun.
Math. Phys. 175 (1996)  427-521; \\
 {\it Geometry of quantum principal 
bundles II},
Rev. Math. Phys. 9 (1997) 531-607.

\bibitem{dgh} L. D\c{a}browski, H. Grosse, P.M. Hajac
\emph{Strong connections and Chern-Connes pairing in the Hopf-Galois theory},
Commun. Math. Phys. 206 (1999) 247-264.

\bibitem{hP} P.M. Hajac,
\emph{Strong connections on quantum principal bundles}, Commun. Math. Phys. 182
(1996) 579-617.

\bibitem{HM} P. M. Hajac, S. Majid \emph{Projective module description of the
   $q$-monopole}, Commun. Math. Phys. 206 (1999) 247-264.

\bibitem{HMS} P. M. Hajac, R. Matthes, W. Szyma\'nski
\emph{A locally trivial quantum Hopf fibration}, 
arXiv:math.QA/0112317; to appear
in Algebra and Representation Theory.


\bibitem{la00} G. Landi,
\emph{\it Deconstructing monopoles and instantons}, Rev. Math. Phys. 12 (2000)
1367-1390.

\bibitem{maj} S. Majid,
\emph{Quantum and braided group Riemannian geometry}, J. Geom. Phys.  30 (1999) 113-146.

\bibitem{ka} C. Kassel, {\it Quantum groups}, Springer 1995.

\bibitem{KS} A. Klimyk, K.  Schm\"udgen,
\emph{Quantum groups and their  representations}, Springer-Verlag 
Berlin Heidelberg
1997.

\bibitem{KT} H. F. Kreimer, M. Takeuchi,
\emph{Hopf algebras  and Galois extensions of an algebra}, Indiana 
Univ. Math. J.
30  (1981) 675-692.

\bibitem{MNW} T. Masuda, Y. Nakagami, J. Watanabe,
\emph{Noncommutative differential geometry on the quantum $SU(2)$. I:An algebraic
viewpoint}, K-Theory 4 (1990) 157-180; \\
\emph{Noncommutative differential geometry on the quantum two sphere of P.Podle{\'s}. I: An algebraic
viewpoint}, K-Theory 5 (1991) 151-175.

\bibitem{HG} S. Montgomery,
\emph{Hopf algebras  and their actions on rings},  AMS 1993.

\bibitem{future} C. Pagani, in preparation.

\bibitem{po} P. Podle{\'s},
\emph{Quantum spheres},  Lett.\ Math.\ Phys.\  14 (1987)  193-202.

\bibitem{frt}  N.Yu. Reshetikhin, L.A. Takhtadzhyan, L.D.  Faddeev,
\emph{Quantization of Lie groups and Lie algebras},  Leningrad Math. 
J.  1 (1990)
193-225.

\bibitem{scP} P. Schauenburg,
\emph{Bi-Galois objects over Taft algebras},  Israel J. Math. 
115 (2000) 101-123.

\bibitem{scsc} P. Schauenburg, H. Schneider,
\emph{Galois type extensions of noncommutative algebras},  in  preparation.

\bibitem{scH} H. Schneider,
\emph{Principal  homogeneous spaces for arbitrary Hopf algebras},  Israel J. of
Math.  72 (1990) 167-195.

\bibitem{si} B. Simon,  {\it Trace ideals and  their applications},  Cambridge
Univ. Press, 1979.


\bibitem{wo} S.L. Woronowicz,
\emph{Twisted $\SU(2)$ group.  An example of a noncommutative differential
calculus},  Publ.\ Res.\  Inst.\ Math.\ Sci.\ 23 (1987)  117-181.


\end{thebibliography}
\end{document}